\documentclass{elsarticle}%,openany

\usepackage{amsfonts,amssymb,amsmath,amsthm}%,naglowki}
\usepackage[dvips]{color}

\def\num{\hspace{-2mm}{\bf }\hspace{2mm}}

\newtheorem{st}{Statement}[section]

\newtheorem{prob}[st]{Problem}

\newtheorem{propo}[st]{Proposition}
\newtheorem{ex}[st]{Example}
\newtheorem{cor}[st]{Corollary}

\newtheorem{thm}[st]{Theorem}

\def\max{{\rm max\,}}
\def\min{{\rm min\,}}

\def\card{{\rm card\,}}

\def\dom{{\rm dom\,}}
\def\indec{{\rm indec\,}}
\def\rng{{\rm rng\,}}
\def\ot{{\rm ot\,}}
\def\prede{{\rm pred\,}}

\def\Proof:{ \vspace{-1.5mm} {\noindent\it Proof.}}
\def\gd{\vspace{0mm}}%-5mm

\def\Box{\rule{1.5mm}{1.5mm}}

\begin{document}

\renewcommand{\thefootnote}{\fnsymbol{footnote}}

 \title{ Cascades, Order and Ultrafilters}
 \author{ Andrzej Starosolski}
 \ead{Andrzej.Starosolski@polsl.pl}
 \address{Politechnika \'Sl\c{a}ska, Wydzia{\l} Matematyki Stosowanej ul. Kaszubska 23, 44-101 Gliwice, Poland}
 \date{\today}

\begin{abstract}

We investigate mutual behavior of cascades, contours of which are
contained in a fixed ultrafilter. Using that relation we prove (ZFC)
that the class of strict $J_{\omega^\omega}$-ultrafilters,
introduced by J. E. Baumgartner in \textit{Ultrafilters on $\omega$}, is empty. We
translate the result to the language of $<_\infty$-sequences under
an ultrafilter, investigated by C. Laflamme in \textit{A few special ordinal ultrafilters}, to show
that if there is an arbitrary long finite $<_\infty$-sequence under
$u$ than $u$ is at least strict $J_{\omega^{\omega+1}}$-
ultrafilter.
\end{abstract}

\maketitle

\gd
\section{Introduction}

%\footnotetext{2010 {\it Mathematics Subject Classification} 03E04; 03E05.

%{\it Key words and phrases}: ordinal ultrafilters; monotone sequential cascades.}

Baumgartner in the article \textit{Ultrafilters on $\omega$} (\cite{Baum}) introduced a notion of
${\cal I}$-ultrafilters: Let $\cal I$ be an ideal on $X$, an ultrafilter (on $\omega$)
is an ${\cal I}$-ultrafilter, if and only if,
for every function $f:\omega \rightarrow X$ there is a set $U\in u $ such that $f[U] \in {\cal I}$.
This kind of ultrafilters was studied by large group of mathematician. We shall mention only the most important
papers in this subject from our point of view: J. Brendle \cite{Brendle}, C. Laflamme
\cite{Laf}, Shelah \cite{Shelah2} ,\cite{Shelah}, B\l aszczyk \cite{Blaszczyk}.  
Among other types of ultrafilters J. E. Baumgartner introduced ordinal
ultrafilters, precisely $\omega_1$ sequence of classes of ultrafilters.
We say that $u$ is $J_\alpha$ ultrafilter (on $\omega$)
if for each function $f:\omega \rightarrow \omega_1$ there is $U \in
u$ such that $\ot(f(U)) <\alpha$, where $\ot(\cdot )$ denotes the ordre type. For additional information about
ordinal ultrafilters a look at \cite{Baum}, \cite{Brendle},
\cite{Star-P-hier} is recommended. In \cite{Baum} J. E. Baumgartner
proved (in Theorems 4.2 and 4.6) that for each successor ordinal
$\alpha < \omega_1$ the class of strict $J_{\omega^{\alpha}}$-ultrafilters
(see below) is nonempty if P-points exist, he also
pointed out that: "In general we do not know, whether, if $\alpha$
is limit, there is a $J_{\omega^\alpha}$-ultrafilter that is no
$J_{\beta}$-ultrafilter, for some $\beta < \omega^\alpha$, even if
CH or MA assumed". Here, such ultrafilters we call strict
$J_{\omega^\alpha}$-ultrafilters, and we partially solve the
problem, showing (ZFC) that the class of strict
$J_{\omega^\omega}$-ultrafilters is empty.

If $u$ is a filter(base) on $A \subset B$, then we identify $u$ with
the filter on $B$ for which $u$ is a filter-base. Let $u$, $v$ be
ultrafilters on $\omega$, recall that $v <_\infty u$ if there is a
function $f:\omega \rightarrow \omega$ such that $f(u)=v$ and $f$ is
not finite-to-one or constant on any set $U \in u$. In \cite{Laf} C.
Laflamme proved (reformulation of Lemma 3.2) that if an ultrafilter $u$ has an
infinite decreasing $<_\infty$- sequence below, then $u$ is at least
strict $J_{\omega^{\omega+1}}$-ultrafilter. He also stated the following

\cite[Open Problem 1]{Laf} What about the corresponding influence of
increasing $<_\infty$-chains below $u$? Given such an ultrafilter
$u$  with an increasing infinite $<_\infty$-sequence $u>_{RK} \ldots
>_\infty u_1 >_\infty u_0$ below, fix maps $g_i$ and $f_i$
witnessing $u >_{RK} u_i$ and $u_{i+1} >_\infty u_i$ respectively.
The problem is really about the possible connections between $g_i$
and $f_i\circ g_{i+1} $ even relative to members of $u$.

\cite[Open Problem 2]{Laf} Can we have an ultrafilter $u$ with
arbitrary long finite $<_\infty$-chains below $u$ without infinite
one? This looks like the most promising way to build a strict
$J_{\omega^\omega}$-ultrafilter.

We find affirmative answer to the first problem and  negative
answer to the second one.

\section{Prelimineries}

In \cite{DolMyn} S. Dolecki and F. Mynard introduced monotone
sequential cascades - special kind of trees - as a tool to describe
topological sequential spaces. Cascades and their contours appeared
to be also a useful tool to investigate certain types of ultrafilters
on $\omega$, namely ordinal ultrafilters and P-hierarchy
(see \cite{Star-P-hier}, \cite{Star-Ord-vs-P}), here we focus on the first of them.

The {\it cascade} is a tree $V$, ordered by "$\sqsubseteq $",
 without infinite branches and with a
least element $\emptyset _V$. A cascade is $\it sequential$ if for
each non-maximal element of $V$ ($v \in V \setminus \max V$) the set
$v^{+V}$ of immediate successors of $v$ (in $V$) is countably
infinite. We write $v^+$ instead of $v^{+W}$ if it is known in which
cascade the successors of $v$ are considered. If $v \in V \setminus
\max V$, then the set $v^+$ (if infinite) may be endowed with an
order of the type $\omega$, and then by $(v_n)_{n \in \omega}$ we
denote the sequence of elements of $v^+$, and by $v_{nW}$ - the
$n$-th element of $v^{+W}$. We say that $v$ is a {\it predecessor}
of $v'$ (we write $v=\prede(v')$)
if $v' \in v^+$.

The {\it rank} of $v \in V$ ($r_V(v)$ or $r(v)$) is defined
inductively as follows: $r(v)=0$ if $v \in \max V$, and otherwise
$r(v)$ is the least ordinal greater than the ranks of all immediate
successors of $v$. The rank $r(V)$ of the cascade $V$ is, by
definition, the rank of $\emptyset_V$. If it is possible to order
all sets $v^+$ (for $v \in V \setminus \max V$)  so that for each $v
\in V \setminus \max V$ the sequence $(r(v_n)_{n<\omega})$ is
non-decreasing (other words if for each 
$v\in V\setminus \emptyset_V$ the set $\{v'\in (\prede(v))^+
: r(v')<\alpha\}$ is finite for each $\alpha<r(v)$),
then the cascade $V$ is {\it monotone}, and we fix
such an order on $V$ without indication. Thus we introduce
lexicographic order $<_{lex}$ on $V$ in the following way:
$v' <_{lex} v''$ if $v'\sqsupset v''$ or if there exist $v$,
$\tilde{v'} \sqsubseteq v'$ and $\tilde{v''}\sqsubseteq v''$
such that $\tilde{v'}\in v^+$ and $\tilde{v''}\in v^+$ and
$\tilde{v'}=v_n$, $\tilde{v''}=v_m$ and $n<m$.

Let $W$ be a cascade, and let $\{V_w: w \in \max W\}$ be a set of
pairwise disjoint cascades such that $V_w \cap W = \emptyset$ for
all $w \in \max W$. Then, the {\it confluence} of cascades $V_w$
with respect to the cascade $W$ (we write $W \looparrowleft V_w$) is
defined as a cascade constructed by the identification of $w \in
\max W$ with $\emptyset_{V_w}$ and according to the following rules:
$\emptyset_W = \emptyset_{W \looparrowleft V_w}$; if $w\in W
\setminus \max W$, then $w^{+ W \looparrowleft V_w} = w^{+W}$; if $w
\in V_{w_0}$ (for a certain $w_0 \in \max W$), then $w^{+ W
\looparrowleft V_w} = w^{+V_{w_0}}$; in each case we also assume
that the order on the set of successors remains unchanged. By $(n)
\looparrowleft V_n$ we denote $W \looparrowleft V_w$ if $W$ is a
sequential cascade of rank 1.

Also we label elements of a cascade $V$ by sequences of naturals of
length $r(V)$ or less, by the function which preserves the lexicographic
order, $v_l$ is a resulting name for an element of $V$, where $l$ is
the mentioned sequence (i.e. $v_{l^\frown n}= (v_l)_n \in v^+$); by
$V_l$ we denote $v_l^\uparrow$ and by $L_{\alpha,V}$ we understand
$\{l\in \omega^{<\omega}: r_V(v_l)=\alpha\}$. Let $v$, $v' \in V$,
we say that $v'$ is a {\it predecessor} of $v$ (in $V$) if $v \in
v'^+$, we write $v'=\prede_V(v)$. For a finite sequence
$l=(n_0, \ldots , n_k)$ of natural numbers
by $l^-$ we denote a sequence $l$ with the last
element removed, i.e. $l^-=(n_0, \ldots , n_{k-1})$; by $l^+$
we denote a set of all sequences $l'$ such that $l'^-=l$.

If $\mathbb{U}= \{{u}_s: s \in S \}$ is a family of filters on $X$
and if $p$ is a filter on $S$, then the {\it contour of $\{ {u}_s
\}$ along $p$} is defined by
$$\int_{p} \mathbb{U} = \int_{p}{u}_s =
\bigcup_{P \in p} \bigcap_{s \in P} {u}_s.$$

Such a construction has been used by many authors (\cite{Fro},
\cite{Gri1}, \cite{Gri2}) and is also known as a sum (or as a limit)
of filters. On the sequential cascade, we consider the finest
topology such that for all but the maximal elements $v$ of $V$, the
co-finite filter on the set $v^{+V}$ converges to $v$. For the
sequential cascade $V$ we define the {\it contour} of $V$ (we write
$\int V$) as the trace on $\max V$ of the neighborhood filter of
$\emptyset _V$ (the {\it trace} of a filter $u$ on a set $A$ is the
family of intersections of elements of $u$ with $A$). Similar
filters were considered in \cite{Kat1}, \cite{Kat2}, \cite{Dagu}.
Let $V$ be a monotone sequential cascade and let $u=\int V$.
Then
the {\it rank $r(u)$} of $u$
%($r(u)$)
is, by definition, the rank of $V$.
%, it
It
was shown in \cite{DolStaWat} that if $\int V= \int W$, then $r(V) =
r(W)$.
%).

Let $S$ be a countable set. A family $\{u_s\}_{s\in S}$ of filters
is referred to as $\it discrete$ if there exists a pairwise disjoint
family $\{U_s\}_{s \in S}$ of sets such that $U_s \in u_s$ for each
$s \in S$. For $v \in V$ we denote by $v^\uparrow$ a subcascade of
$V$ built by $v$ and all successors of $v$. If $U \subset \max V$
and $U \in \int V$, then by $U^{\downarrow V}$ we denote the biggest
(in the set-theoretical order) monotone sequential subcascade of
cascade $V$ built of some $v \in V$ such that $U \cap \max v
^\uparrow \not= \emptyset$.  We write $v^{\uparrow }$ and
$U^\downarrow$ instead of $v^{\uparrow V}$ and $U^{\downarrow V}$ if
we know in which cascade the subcascade is considered. The reader may find more
information about monotone sequential cascades and their contours in \cite{Dol-multi}, \cite{DolMyn},
\cite{DolStaWat}, \cite{Star-ff}, \cite{Star-Ord-vs-P},
\cite{Star-P-hier}.

In the remainder of this paper each filter is considered to be on
$\omega$, unless indicated otherwise.

\section{Existence of ordinal ultrafilters}

For a monotone sequential cascade $V$ by $f_V$ we denote an {\it
lexicographic order respecting function} $\max V \rightarrow \omega_1$, i.e., such
a function that $v'<_{lex}v''$ iff $f_V(v')< f_V(v'')$ for each
$v'$, $v'' \in \max V$. If $f:\omega \rightarrow \omega_1$ and $f=f_V$ for
some monotone sequential cascade $V$ then we say that $V$
corresponds with an order of $f$.

Let $V$ and $W$ be monotone sequential cascades such that $\max V
\supset \max W$. We say that $W$ {\it increases the order of} $V$
(we write) $W\Rrightarrow V$ if $\ot(f_W(U)) \geq
\indec(\ot(f_V(U)))$ for each $U \subset \max W$, where
$\indec(\alpha)$ is the biggest indecomposable ordinal less then, or
equal to $\alpha$; by Cantor normal form theorem such a number
exists and is defined uniquely. Clearly this relation is idempotent
and transitive. Although  relation of increasing of order says
 that one cascade is somehow bigger then another,
this relation is quite independent with the containment of contours.

\begin{ex}($\int T \supset \int V \not\Rightarrow  T\Rrightarrow V$)
Let $(V_n)_{n\in \omega}$ be a sequence of pairwise disjoint monotone sequential cascades of rank 2.
For each $n<\omega$ choose $v_n$ - an arbitrary element of $\max V_n$. Let $V'_n=V_n \setminus \{v_n\}$
 and let $V'_\omega$
be an arbitrary monotone sequential cascade of rank 2 such that $\max V_\omega = \bigcup_{n<\omega}\{v_n\}$. Now put
$T = (n) \looparrowleft_{n < \omega}V_n$ and $V = (n) \looparrowleft_{n \leq \omega}V'_n$.
\end{ex}

\begin{ex} ($T\Rrightarrow V \not\Rightarrow \int T \supset \int V$)
Let $(V_n)_{n\in \omega}$ be a sequence of pairwise disjoint monotone sequential cascades of rank 1.
For each $n<\omega$ choose $v_n$ - an arbitrary element of $\max V_n$. Let $(B_n)_{n<\omega}$ be a partition of
$\bigcup_{n<\omega}\{v_n\}$ into infinite sets. Let $V'_n$ be a monotone sequential cascade of rank 1,
such that $\max V'_n = (\max V_n \setminus \{v_n\})\cup B_n$. Put $T= (n) \looparrowleft V'_n$ and
$V= (n) \looparrowleft V_n$.
\end{ex}

Let $V$ and $W$ be monotone sequential cascades, let
$f: V\rightarrow W$ be a finite-to-one, $\sqsubseteq_V$ order
preserving surjection such that $F_{\mid \max V \cup \emptyset_V} =
id_{\mid \max V \cup \emptyset_V}$ and $v \in f^{-1}(v)$ for each $V
\in V$. Then, it is easy to see, that $V\Rrightarrow f(V)$ and $f(V)
\Rrightarrow V$, we call this property locally finite partition
property (LFPP)

Let $u$, $p$ be filters on $\omega$, then $u\vee p$ we define as
$\{X\in \omega:$ there exist $U\in u$ and $P\in p$ such that $U\cap
P \subset X\}$.

Let $u$ be an ultrafilter and let $V$, $W$ be monotone sequential
cascades such that $\int V \subset u$ and $\int W \subset u$. Then
we say that {\it rank $\alpha$ in cascade $V$ agree with rank
$\beta$ in cascade $W$ with respect to the ultrafilter $u$}  if for
any choice of $\tilde{V}_{p,s} \in \int V_p$ and $\tilde{W}_{p,s}
\in \int W_s$ there is: $\bigcup_{(p,s)\in L_{\alpha,V}\times
L_{\beta,W}} ( \tilde{V}_{p,s} \cap \tilde{W}_{p,s}) \in u$; this
relation is denoted by $\alpha_VE_u\beta_W$.

\begin{propo}\num
Let $u$ be an ultrafilter and let $V$, $W$ be monotone sequential
cascades such that $\int V \subset u$ and $\int W \subset u$. Then
$1_VE_u1_W$ and $r(V)_VE_ur(W)_W$.
\end{propo}

\Proof: First suppose that $r(V)=1$ or $r(W)=1$, say $r(V)=1$.
Clearly $\card(L_{1,V})=1$. For each $(p,s)\in L_{1,V} \times
L_{1,W}$ take any $\tilde{V}_{p,s} \in \int V_p$ and any
$\tilde{W}_{p,s} \in \int W_s$. Notice that since $\int V$ is a
co-finite filter on $\max V$ thus for each $(p,s)\in L_{1,V} \times
L_{1,W}$ the set $(\tilde{W}_{p,s} \cap \max V) \setminus
(\tilde{W}_{p,s} \cap \tilde{V}_{p,s})$ is finite. Therefore there
exist $\hat{W}_{p,s}\subset \tilde{W}_{p,s}$, $\hat{W}_{p,s}\in \int
W_{p,s}$ such that $\hat{W}_{p,s}\cap \tilde{V}_{p,s} =
\hat{W}_{p,s}\cap \max V$. Thus $\bigcup_{(p,s)\in L_{1,V} \times
L_{1,W}}(\tilde{W}_{p,s} \cap \tilde{V}_{p,s}) \supset
\bigcup_{(p,s)\in L_{1,V} \times L_{1,W}}(\hat{W}_{p,s} \cap
\tilde{V}_{p,s}) = \bigcup_{(p,s)\in L_{1,V} \times
L_{1,W}}(\hat{W}_{p,s}) \cap \max V$. Since $\bigcup_{(p,s)\in
L_{1,V} \times L_{1,W}}(\hat{W}_{p,s})\in \int W \subset u$ and
$\max V \in \int V \subset u$ thus $\bigcup_{(p,s)\in L_{1,V} \times
L_{1,W}}(\tilde{W}_{p,s} \cap \tilde{V}_{p,s}) \in u$.

Before we deal with case $r(V)\geq 2$ and $r(W)\geq 2$ we state the
following claim: in assumption of this Proposition, if a set $U$ is
such that for each $(p,s)\in L_{1,V} \times L_{1,W}$ the
intersection $U\cap \max V_p \cap \max W_s$ is finite, then $U
\not\in u$. Let $f:\omega \rightarrow L_{1,V}$, $h:\omega
\rightarrow L_{1,W}$ be bijections, and let $G(i,j) = U\cap \max
V_{f(i)} \cap \max W_{h(j)}$ for $i,j \in \omega$. Let
$\Delta^{\geq} = \{(i,j)\subset \omega\times\omega: i\leq j\}$,
$\Delta^{\leq} = \{(i,j)\subset \omega\times\omega: i\geq j\}$.
Since $u$ is an ultrafilter thus either $G(\Delta^\geq) \in u$ or
$G(\Delta^\leq) \in u$. But $G(\Delta^\geq )$  is finite on each
$\max V_p$ for $p\in L_{1,V}$ and so $(G(\Delta^\geq ))^c \in \int
V$ therefore $G(\Delta^\geq ) \not\in u$. Also $G(\Delta^\leq )$  is
finite on each $\max W_s$ for $s\in L_{1,W}$ and so $(G(\Delta^\leq
))^c \in \int W$ therefore $G(\Delta^\leq ) \not\in u$.

Now let $r(V)\geq 2$ and $r(W)\geq 2$ and suppose on the contrary
that $K=\bigcup_{(p,s)\in L_{1,V} \times L_{1,W}}(\tilde{W}_{p,s}
\cap \tilde{V}_{p,s}) \not\in u$ for some choice of $\tilde{V}_{p,s}
\in \int V_p$ and $\tilde{W}_{p,s} \in \int W_s$. Thus $K^c\in u$.
Put $R^\#= \{{(p,s)\in L_{1,V} \times L_{1,W}}; \int V_p \# \int W_s
\}$ and $R^{\not\#}= \{{(p,s)\in L_{1,V} \times L_{1,W}} ; \int V_p
\not\# \int W_s \}$. Define $K_1=\bigcup_{ (p,s)\in R^\#} ((\max
V_p\cap \max W_s)\setminus (\tilde{V}_{p,s} \cap \tilde{W}_{p,s}))$
and $K_2=\bigcup_{(p,s)\in R^{\not\#}}(\max V_p \cap \max W_s)$ and
notice that $K^c \subset K_1 \cup K_2$. Thus either $K_1 \in u$, or
$K_2\in u$. Since traces of $\int W_s$ and of $\int V_p$ on $\max
V_p\cap \max W_s$ are co-finite filters (on $\max V_p\cap \max W_s$)
thus $K_1 \cap \max V_p\cap W_s = (\max V_p \cap \max W_s)\setminus
(\tilde{V}_{p,s} \cap \tilde{W}_{p,s})$ is finite on each $(p,s)\in
R^\#$ and empty on each $(p,s) \in R^{\not\#}$. By similar reasoning
$K_2$ is finite on each $(p,s)\in L_{1,V} \times L_{1,W}$. Therefore
by claim above $K_1 \not\in u $ and $K_2 \not\in u$ - contradiction.
Second statement of Proposition 3.3 is clear. $\Box$

\vspace{3mm}
In the above Proposition 3.3 the inverse of the
implication does not hold.

\begin{ex} Let $(A_\alpha)_{\alpha \leq 2\omega}$ be
a partition of $\omega$ into infinite sets. Let $V_\alpha$ be a
monotone sequential cascade of rank 1 such that $\max V_\alpha =
A_\alpha$. Let $W= (\alpha)\looparrowleft_{\alpha \leq \omega}
V_\alpha$, $V= (\alpha) \looparrowleft_{\omega \leq \alpha
<2\omega}V_\alpha$, and let $u$ be any free ultrafilter containing
$A_\omega$. Clearly $\int V \not\subset u$, $\int W \not\subset u$
but $1_VE_u1_W$.
\end{ex}

Although the following Theorem 3.5 is stated using the
"$\Rrightarrow$" relation, it is worth to look at the proof of it as
on the description of possible relations of cascades whose contours
are contained in the same ultrafilter, and as a description of
operation which leads from such cascades to others whose contours
are also contained in the same ultrafilter.

\begin{thm}\num
Let $u$ be an ultrafilter and let $V$, $W$ be monotone sequential
cascades of finite ranks such that $\int V \subset u$ and $\int W
\subset u$. Then $n_VE_um_W$ implies the existence of a monotone
sequential cascade $T$ of rank \linebreak $\max
\{r(V),r(W)\} \leq r(T) \leq r(V)+r(W)$
and such that $\int T \subset u$, $T \Rrightarrow
V$, $T \Rrightarrow W$, $\int V \subset \int T$ and $\int W \subset \int T$.
\end{thm}

\Proof: Before we start the proof, let us make the following remarks:
in this theorem we claim (in place of $\int T \subset u$) that
(under the same assumption and notation) $\int T \subset (\int V
\vee \int W)$, and this formulation also will be used in the proof;
cascade $T$ build in this proof has ranks not less then $\max \{n,m\}+ \max \{r(V)-n, r(W)-m\}$
and not greater then $r(V)+r(W)-1$ and this is inductively used in the proof.
Without loss of generality, we may assume that each branch in $V$
has length $r(V)$ and each branch in $W$ has length $r(W)$, and that $r(V)\leq r(W)$.

We proceed by induction by $r(W)$, and for each $r(W)$ by sub-induction by $r(V)$.
First step of induction and of sub-inductions is $r(V)=1$ and then we take
$T = W^{\downarrow \max V} $ which clearly
fulfills the claim. Assume that the claim is proved for all cascades $\breve{V}$, $\breve{W}$
which behave like in assumptions and such that $r(\breve{V})\leq r(\breve{W})$
and  either $r(\breve{W}) < r({W})$ or else ($r(\breve{W}) = r({W})$ and
$r(\breve{V})< r({V})$).

We consider 3 cases

1) $n<r(V)$ and $m<r(W)$;

2 ($n=r(V)$ and $m<r(W)$) or ($n<r(V)$ and $m=r(W)$);

3) $n=r(V)$ and $m=r(W)$.

Case 1) 
Let  $R^\#=\{ (l,s) \in L_{n,V} \times L_{m,W}: \int V_l \#
\int W_s\}$. Notice that exactly one of the following 3 subcases
holds:

$1.1$) There is $K_{F-F} \subset R^\#$ such that
$\card(K_{F-F}(l))<\omega$ and \linebreak $\card(K^{-1}_{F-F}(s))<\omega$ for
each $l \in \dom K_{F-F}$, $s \in \rng K_{F-F}$ and \linebreak
$\bigcup_{(l,s)\in K_{F-F}}(\tilde{V}_{l,s} \cap \tilde{W}_{l,s})
\in u$ for each choice of $\tilde{V}_{l,s} \in \int V_{l,s}$ and of
$\tilde{W}_{l,s} \in \int W_{l,s}$;

$1.2$) $\sim 1.1$ and there is $K_{\infty-F}\subset R^\#$ such that
$\card(K_{\infty-F}^{-1}(s))<\omega$ for each $s\in \rng
K_{\infty-n}$ and $\bigcup_{(l,s)\in K_{\infty-F}}(\tilde{V}_{l,s}
\cap \tilde{W}_{l,s}) \in u$ for each choice of $\tilde{V}_{l,s} \in
\int V_{l,s}$ and of $\tilde{W}_{l,s} \in \int W_{l,s}$;

$1.3$) $\sim 1.1$ and there is $K_{F-\infty}\subset R^\#$ such that
$\card(K_{F-\infty}(l))<\omega$ for each $l\in \dom K_{F-\infty}$
and $\bigcup_{(l,s)\in K_{F-\infty}}(\tilde{V}_{l,s} \cap
\tilde{W}_{l,s}) \in u$ for each choice of $\tilde{V}_{l,s} \in \int
V_{l,s}$ and of $\tilde{W}_{l,s} \in \int W_{l,s}$.

Let  $\Delta^\leq =\{(i,j):i\geq j; i,j \in \omega\}$ and
$\Delta^\geq =\{(i,j):i\leq j; i,j \in \omega\}$. Let $p:\omega
\rightarrow L_{n,V}$ and $q: \omega \rightarrow L_{m,W}$ be
bijections. For $X \subset \omega \times \omega$ define $G(X)=
\bigcup_{(i,j)\in X} (\max V_{p(i)} \cap \max W_{q(j)})$ and put
$(p,q)(x_1,x_2) = (p(x_1),q(x_2))$. Since $G(\Delta^\geq) \cup
G(\Delta^\leq) \in u$ thus either $G(\Delta^\geq)\in u$ (case $1.1$
or $1.2$) or $G^\geq(\Delta) \in u$ (case $1.1$ or $1.3$). For $A
\subset R$ we define also $H(A)=\bigcup_{(l,s)\in A}(\max V_l \cap
\max W_s)$.

Since $u$ is an ultrafilter thus without loss of generality (by
LFPP,  for case 1.1 used twice, its property of increasing order and
transitivity of "$\Rrightarrow$" relation) exactly one of the
following subcases holds.

$1.1'$) There is $K_{1-1} \subset R^\#$ such that $
\card(K_{1-1}(l))=1$ and $\card(K^{-1}_{1-1}(s))=1$ for each
$(l,s)\in K_{1-1}$and $\bigcup_{(l,s)\in K_{1-1}}(\tilde{V}_{l,s}
\cap \tilde{W}_{l,s}) \in u$ for each choice of $\tilde{V}_{l,s} \in
\int V_{l,s}$ and of $\tilde{W}_{l,s} \in \int W_{l,s}$;

$1.2'$) $\sim 1.1$ and there is $K_{\infty-1}\subset R^\#$ such
that $\card(K_{\infty-1}^{-1}(s))=1$ and
$\card(K_{\infty-1}(l))=\omega$ for each $(l,s) \in K_{\infty-1}$
and $\bigcup_{(l,s)\in K_{\infty-1}}(\tilde{V}_{l,s} \cap
\tilde{W}_{l,s}) \in u$ for each choice of $\tilde{V}_{l,s} \in \int
V_{l,s}$ and of $\tilde{W}_{l,s} \in \int W_{l,s}$;

$1.3'$) $\sim 1.1$ and there is $K_{1-\infty}\subset R^\#$ such that
$\card(K_{1-\infty}(l))=1$ and $\card(K^{-1}_{1-\infty}(s))=\omega$
for each $(l,s)\in K_{1-\infty}$ and $\bigcup_{(l,s)\in
K_{1-\infty}}(\tilde{V}_{l,s} \cap \tilde{W}_{l,s}) \in u$ for each
choice of $\tilde{V}_{l,s} \in \int V_{l,s}$ and of $\tilde{W}_{l,s}
\in \int W_{l,s}$.

 Subcase $1.1'$) Without loss of generality, we may assume that
$\max V = \max W = H(K_{1-1})$. Define a series of sets:
$R= L_{n,V} \times L_{m,W}$, \linebreak
$R(A) = \{ (l,s) \in R: \card(\max V_l \cap \max W_s \cap A) =
\omega\}$ for $A \subset \omega$, $\hat{u} = \{ R(U): U \in u\}$,
$\hat{V}=\{v\in V: r(v) >n\}\cup R(\max V)$, $\hat{W}=\{w\in W:
r(w)>m\}\cup R(\max W)$. On $\hat{V}$ we define order
$\sqsubseteq_{\hat{V}}$ by: if $v_1$, $v_2 \in V \cap \hat{V}$,
$(l,s) \in R$ then: $v_1\sqsubseteq_V v_2$ iff
$v_1\sqsubseteq_{\hat{V}} v_2$; $v_1 \sqsubseteq_{\hat{V}} (l,s)$
iff $\card(\max V_l \cap \max W_s \cap \max v_1^{\uparrow
V})=\omega$.  On $\hat{W}$ we intrtoduce order in the analogical
way. Notice that $\hat{u}$ is an ultrafilter on $R$ and that
$\hat{V}$ and $\hat{W }$ are monotone sequential cascades (on $R$)
and that $r(\hat{ V})=r(V)-n$, $r(\hat{W})= r(W)-m$ and that $\int
\hat{V} \subset \hat{u}$ and $\int \hat{{W}} \subset \hat{u}$.

By inductive assumption there is $\hat{T}$ monotone sequential
cascade (on $R$) of rank
$\max\{r(V)-n, r(W)-m\} \leq r(\hat{T}) \leq r(V)+r(W)- n - m - 1$ and
such that $\int \hat{T} \subset \hat{u}$ and $\hat{T} \Rrightarrow
\hat{V}$, and $\hat{T} \Rrightarrow \hat{W}$, also by inductive
assumption, for each $(l,s) \in R^\# $ there is a monotone sequential
cascade $T_{l,s}$ of rank $\max\{n,m\} \leq  r(T_{l,s}) \leq n+m-1$
such that $T_{l,s} \Rrightarrow V_{l}$, 
$T_{l,s} \Rrightarrow W_{s}$, $\int T_{l,s}
\subset \int V_l \vee \int W_s$. Define $\bar{T} = \hat{T}
\looparrowleft_{\{(l,s); \max V_l \cap \max W_s \in \max \hat{T}\}}
T_{l,s}$. Take any $A \in \int \bar{T}$, thus there exist $\hat{A} \in
\int\hat{T}$ and $A_{l,s} \in \int T_{l,s}$ such that $A =
\bigcup_{(l,s)\in \hat{A}}A_{(l,s)}$. Since $A_{(l,s)} \supset
\tilde{V}_l \cap \tilde{W}_s$ for some $\tilde{V}_l \in \int V_l$
and $\tilde{W}_s \in \int W_s$ so $A \supset \bigcup_{(l,s)\in
\hat{A}}(\tilde{V}_{(l,s)}\cap \tilde{W}_{(l,s)})$ $=
\bigcup_{(l,s)\in R}(\tilde{V}_{(l,s)}\cap \tilde{W}_{(l,s)})\cap
\bigcup_{(l,s) \in \hat{A}} (\max V_l \cap \max W_s)$. Since
$\bigcup_{(l,s) \in \hat{A}} (\max V_l \cap \max W_s) \in u$ thus
$\bigcup_{{(l,s)}\in \hat{A}}(\tilde{V}_{(l,s)}\cap
\tilde{W}_{(l,s)})$ $= \bigcup_{(l,s)\in R}(\tilde{V}_{(l,s)}\cap
\tilde{W}_{(l,s)})\cap \bigcup_{(l,s) \in \hat{A}} (\max V_l \cap
\max W_s) \in u$ and so $A \in u$ and so $\int \bar{T} \subset u$.
Consider sets $U_i=\bigcup_{\{(l,s)\in R^\#: r(T_{(l,s)})=i\} }\max
T_{(l,s)}$. By inductive assumption - upper limitation of ranks,
only finite number of these sets are nonempty, and since
$\bigcup_{i<\omega}U_i \in u$ thus $U_{i_0} \in u$ for some $i_0$.
Let $T = \bar{T}^{\downarrow U_{i_0}}$. Clearly $\int T \in u$.
Calculation of the rank of $T$ follows easily. Take any $P\in \max
T$, without loss of generality, we may assume that
$\ot(f_V(P))=\omega^b$ for some $b\leq r(V)$. Split $R$ into
following sets $R_a =  \{(l,s)\in R:\omega^{a-1} \leq \ot(f_V(P\cap \max
V_l \cap \max W_s))< \omega^{a}\}$ for $a \in \{1, \ldots,
b+1\}$. For some $a$, say $a_0$, we have  $\ot(f_V(\bigcup_{(l,s)\in
R_{a_0}} \max V_l \cap \max W_s \cap P) = \omega^b$, thus
$\ot(f_{\tilde{V}}(R_{a_0})\geq \omega^{b-a}$. Therefore
$\ot(f_{\hat{T}}(R_{a_0}))\geq \omega^{b-a}$ and so $\ot(f_T(P))\geq
\omega^{b}$, and so $T\Rrightarrow V$. Proof that $T \Rrightarrow W$
is analogical.

Subcase $1.2'$)  Without loss of generality, we may assume that $\max
V = \max W =  H(K_{\infty-1})$. Consider cascade $V'$ - such a
modification of cascade $V$ that in the place of the cascade $V_l$,
for each $l
\in \dom K_{\infty-1}$ we put a following cascade: %\linebreak
$(s) \looparrowleft_{\{(s): (l,s)
\in K_{\infty-1}\}} V_{l}^{\downarrow  \max W_s \cap
\bigcup_{k>q^{-1}(s)}\max V_{l^\frown k}}$.    
%Consider also sets $K^{\infty,1}=\{(l,s)\in K_{\infty-1}: \card(
%\{(l',s)\in K_{\infty-1}:l'\in (l^-)^+\})=\omega\}$,
%$K^{\infty,j}= \{(l,s) \in \omega^{<\omega} \times \omega^{<\omega}:
%\exists _{(l',s):l'\in l^+, s\in \rng K_{\infty-1}} (l',s)\in K^{\infty,j-1}\Longrightarrow
%\exists^\infty _{(l',s):l'\in l^+, s\in \rng K_{\infty-1}} (l',s)\in K^{\infty,j-1}\} $, for $1<j<r(V)-n$.
%Now take $K^\infty=\{(l,s)\in K_{\infty-1}: \left[ (l,s)\in K^{\infty,j}\Longrightarrow 
%(l^-,s)\in K^{\infty,j+1} \right]$ for all $j<r(V)-n-1\}$.
Notice that
$H(K^\infty)\in \int V$ and so $H(K^\infty) \in u$
so without loss of generality we may assume that 
$K^\infty= K_{\infty-1}$ and so $V'=V'^{\downarrow H(K^\infty)}$.
Notice that
$V'$ is a monotone sequential cascade of rank $r(V')= r(V)+1$ and that $\int V' \in u$.
Calculation of the rank is straightforward, so take $P\in \int V'$
and for each $l\in \dom K_{\infty-1}$ label elements
of the set $\{s: (l,s) \in K_{\infty-1}\}$ by natural numbers by preserving
lexicographic order bijection, $s_n$ is a resulting name.
If $P \in \int V'$ then there exists $\hat{P}\in \int(V\mid_{v\in V: r(v) \leq n})$,
also for each $l: v_l \in \hat{P}$ there exists a co-finite subset
$A_l$ of $\omega$ that for each $(l,s_n)$, such that $v_l\in \hat{P}, n\in A_l$,
there is a set $P_{l^\frown s_n}\in \int  V_{l}^{\downarrow  \max W_s \cap
\bigcup_{k>q^{-1}(s)}\max V_{l^\frown k}}$, such that
$P=\bigcup_{l:v_l\in\hat{P}} \bigcup_{n\in A_l} P_{l^\frown s_n}$.
Since for each pair $(l,s_n)$ there exist sets $\tilde{V}_{l,s_n} \in \int V_l$ and 
$\tilde{W}_{l,s_n} \in \int W_{s_n}$ such that
$P_{l^\frown s_n} \supset \tilde{V}_{l,s_n} \cap \tilde{W}_{l,s_n}$ thus 
$P \supset \bigcup_{l:v_l\in\hat{P}} \bigcup_{n\in A_l} (\tilde{V}_{l,s_n} \cap \tilde{W}_{l,s_n})$. Clearly 
$\bigcup_{l\in \dom K_{\infty_1}: v_l \not\in \hat{P}} \max V_l  \not\in u$ so 
$\bigcup_{l\in \dom K_{\infty_1}: v_l \in \hat{P}}\max V_l \in u$.
On the other hand (by assumption $\sim 1.1$)
$\bigcup_{l\in \dom K_{\infty-1}}\bigcup_{l^\frown s_n: n \not\in A_l}
(\tilde{V}_{l,s_n} \cap \tilde{W}_{l,s_n}) \not\in u$, where
$\tilde{W}_{l,s_n}=\max W_{s_n}$,  $\tilde{V}_{l,s_n}=\max V_{l}$ and
$A_l=\omega$ for $(l,s_n)\in K_{\infty-1}$
such that $v_l \not\in \hat{P}$. Thus
$\bigcup_{l\in \dom K_{\infty-1}}\bigcup_{l^\frown s_n: n \in A_l}
(\tilde{V}_{l,s_n} \cap \tilde{W}_{l,s_n}) \in u$, therefore since
$P\supset \bigcup_{l\in \dom K_{\infty_1}: v_l \in \hat{P}}
\max V_l \cap \bigcup_{l\in \dom K_{\infty-1}}
\bigcup_{l^\frown s_n: n \in A_l}
(\tilde{V}_{l,s_n} \cap \tilde{W}_{l,s_n})$ 
we have $P\in u$ and so $\int V' \in u$.

 We will show that also $V'\Rrightarrow V$ holds.
Take any $A \subset \max V'$ and notice that it suffices to
prove $\ot(f_{V'}(A\cap \max V_l)) \geq \indec (\ot(f_{V}(A\cap \max
V_l))$ for such $V_l$ that $r_V(v_l)=n$. So we fix such $l$ and
consider $A \cap \max V_l$ assuming, without loss of generality,
that $\ot(f_{V}(A\cap \max V_l) =\omega^c$ for some $c\leq n$.
Consider a following sequence of sets $(\max V_{l^\frown k}\cap
A)_{k<\omega}$, there is $k_0$ that $\ot(f_V(A\cap V'_{l^\frown
k_0}))= \omega^c$ or there is $O$ - infinite subset of $\omega$ such
that  $\ot(f_V(A\cap V_{l^\frown k}))= \omega^{c-1}$ for each $k\in
O$. Notice that each $V_{l^\frown k}$ is split, during the
construction of $V'$, into finitely many pieces by sets $\max W_s
\cap \bigcup_{k>q^{-1}(s)}\max V_{l^\frown k}$. So there is $s_0$
such that $\ot(f_{V}(A\cap \max V_{l^\frown k}))=$ $\ot (f_{V}(A\cap
\max V_{l^\frown k} \cap \max W_{s_0}))$ $= \ot (f_{V'}(A\cap \max
V_{l^\frown k} \cap \max W_{s_0}))$. Therefore either $\ot(f_{V'}(A
\cap \max V_{l^\frown k})) = \omega^c$ for some $k<\omega$, or
$\ot(f_{V'}(A \cap \max V_{l \frown k})) \geq \omega^{c-1}$ for
infinite number of $k$'s. Thus $\ot(f_{V'}(A \cap \max V_{l})) \geq
\omega^c$ and so $V' \Rrightarrow V$.

We notice that for cascades $V'$ and $W$ conditions described as 1.1 hold.
Now we proceed like in subcasce 1.1'.
Define a series of sets:
$R'= L_{1,V'} \times L_{1,W}$, $R'^\#=\{(l,s) \in R': \int V'_l\# \int W_s$ %\linebreak
$R'(A) = \{ (l,s) \in R: \card(\max V'_l \cap \max W_s \cap A) =
\omega\}$ for $A \subset \omega$, $\hat{u} = \{ R(U): U \in u\}$,
$\hat{V'}=\{v\in V': r(v) >1\}\cup R(\max V')$, $\hat{W}=\{w\in W:
r(w)>1\}\cup R(\max W)$. Observe that
$(l,s_n) \in R \Longleftrightarrow (l^\frown n, s) \in R'$ and 
$(l,s_n) \in R^\# \Longleftrightarrow (l^\frown n, s) \in R'^\#$. On $\hat{V'}$ we define order
$\sqsubseteq_{\hat{V'}}$ by: if $v_1$, $v_2 \in V' \cap \hat{V'}$,
$(l,s) \in R'$ then: $v_1\sqsubseteq_{V'} v_2$ iff
$v_1\sqsubseteq_{\hat{V'}} v_2$; $v_1 \sqsubseteq_{\hat{V'}} (l,s)$
iff $\card(\max V'_l \cap \max W_s \cap \max v_1^{\uparrow
V'})=\omega$.  On $\hat{W}$ we intrtoduce order in the analogical
way. Notice that $\hat{u}$ is an ultrafilter on $R'$ and that
$\hat{V'}$ and $\hat{W }$ are monotone sequential cascades (on $R'$)
and that $r(\hat{ V'})=r(V')-1=r(V)$, $r(\hat{W})= r(W)-m$ and that $\int
\hat{V'} \subset \hat{u}$ and $\int \hat{{W}} \subset \hat{u}$.

By inductive (or sub-inductive) assumption (for $V'$, $W$ and $u$)
there is $\hat{T}$ monotone sequential
cascade on $R'$ of rank
$\max\{r(V), r(W)-1\} \leq r(\hat{T}) \leq r(V)+r(W)- 2$ and
such that $\int \hat{T} \subset \hat{u}$ and $\hat{T} \Rrightarrow
\hat{V'}$, and $\hat{T} \Rrightarrow \hat{W}$, also by inductive
assumption, for each $(l,s) \in R'^\# $ there is a monotone sequential
cascade $T_{l,s}$ of rank $ r(T_{l,s}) =1 $
such that $T_{l,s} \Rrightarrow V'_{l}$, 
$T_{l,s} \Rrightarrow W_{s}$, $\int T_{l,s}
\subset \int V'_l \vee \int W_s$. Define $\bar{T} = \hat{T}
\looparrowleft_{\{(l,s); \max V'_l \cap \max W_s \in \max \hat{T}\}}
T_{l,s}$. Take any $A \in \int \bar{T}$, thus there exist $\hat{A} \in
\int\hat{T}$ and $A_{l,s} \in \int T_{l,s}$ such that $A =
\bigcup_{(l,s)\in \hat{A}}A_{(l,s)}$. Since $A_{(l,s)} \supset
\tilde{V}'_l \cap \tilde{W}_s$ for some $\tilde{V}'_l \in \int V'_l$
and $\tilde{W}_s \in \int W_s$ so $A \supset \bigcup_{(l,s)\in
\hat{A}}(\tilde{V}'_{(l,s)}\cap \tilde{W}_{(l,s)})$ $=
\bigcup_{(l,s)\in R}(\tilde{V}'_{(l,s)}\cap \tilde{W}_{(l,s)})\cap
\bigcup_{(l,s) \in \hat{A}} (\max V'_l \cap \max W_s)$. Since
$\bigcup_{(l,s) \in \hat{A}} (\max V'_l \cap \max W_s) \in u$ thus
$\bigcup_{{(l,s)}\in \hat{A}}(\tilde{V}'_{(l,s)}\cap
\tilde{W}_{(l,s)})$ $= \bigcup_{(l,s)\in R'}(\tilde{V}_{(l,s)}\cap
\tilde{W}_{(l,s)})\cap $ \linebreak $\bigcup_{(l,s) \in \hat{A}} (\max V'_l \cap
\max W_s) \in u$ and so $A \in u$ and so $\int \bar{T} \subset u$.
Consider sets \linebreak $U_i=\bigcup_{\{(l,s)\in R^\#: r(T_{(l,s)})=i\} }\max
T_{(l,s)}$. By inductive assumption - upper limitation of ranks,
only finite number of these sets are nonempty, and since
$\bigcup_{i<\omega}U_i \in u$ thus $U_{i_0} \in u$ for some $i_0$.
Let $T = \bar{T}^{\downarrow U_{i_0}}$. Clearly $\int T \in u$.
Calculation of the rank of $T$ follows easily. Take any $P\in \max
T$, without loss of generality, we may assume that
$\ot(f_{V'}(P))=\omega^b$ for some $b\leq r(V)$. Split $R'$ into
following sets $R'_a =  \{(l,s)\in R':\omega^{a-1} \leq \ot(f_{V'}(P\cap \max
V'_l \cap \max W_s))< \omega^{a}\}$ for $a \in \{1, \ldots,
b+1\}$. For some $a$, say $a_0$, we have  $\ot(f_{V'}(\bigcup_{(l,s)\in
R'_{a_0}} \max V'_l \cap \max W_s \cap P) = \omega^b$, thus
$\ot(f_{\tilde{V'}}(R'_{a_0})\geq \omega^{b-a}$. Therefore
$\ot(f_{\hat{T}}(R'_{a_0}))\geq \omega^{b-a}$ and so $\ot(f_T(P))\geq
\omega^{b}$, and so $T\Rrightarrow V'$ and $T\Rrightarrow V$ by transitivity of $\Rrightarrow$ relation.
Proof that $T \Rrightarrow W$
is analogical.

Subcase $1.3'$) Proof is analogical to 1.2'.

Case 2) In both subcases proof is an easier version of proof in case 1 (sub-cases 1.2 and 1.3).

Case 3) a) Case $r(V)=r(W)=1$ was done at the beginning of the proof;

b) If $\min \{r(V),r(W)\}=1$ and $\max \{r(V),r(W)\}>1$ then
$1_VE_u1_W$ by Proposition 3.3, and by already proved part 2 of the
proof, the required cascade $T$ exists;

c) If $\min \{r(V),r(W)\}>1$ then $1_VE_u1_W$ by Proposition 3.3,
and by already proved case 1, the required cascade $T$
exists.

Inclusions of contours is straightforward by monotonicity of contour
operation with respect to the confluence.
 $\Box$
 \\

By Proposition 3.3 and Theorem 3.5 we have
\begin{cor}\num
Let $u$ be an ultrafilter and let $V$, $W$ be monotone sequential
cascades of finite ranks. If $\int V \subset u$ and $\int W \subset u$ then there is
a monotone sequential cascade $T$ of finite rank not less then $\max\{r(V), r(W)\}$ and
such that $\int T \subset u$ and $T\Rrightarrow V$.

\end{cor}

\begin{propo}\num\cite[Proposition 3.3 redefined in virtue of it's proof]{Star-P-hier}
Let $V$ be a monotone sequential cascade of rank $\alpha$. If $u$ is
such an ultrafilter that $\int V \subset u$ then  $\ot(f_V(U)) \geq
\omega^\alpha$ for all $U \in u$.
\end{propo}

\begin{propo}\num\cite[Proposition 3.6]{Star-P-hier}
Let $\alpha $ be a countable indecomposable ordinal, let $n<\omega$
and let $u$ be an ultrafilter. If there is a function $f:\omega
\rightarrow \omega_1$ such that $\ot(f(U)) \geq \omega^{\alpha+n}$
for each $U \in u$ and for each $g:\omega \rightarrow \omega_1$
there is $U_g\in u$ such that $\ot(g(U_g)) < \omega^{\alpha +
\omega}$, then there exists a monotone sequential cascade $V$ of
rank $n$ such that $\int V \subset u$.
\end{propo}

\begin{thm}\num
(ZFC) The class of strict $J_{\omega^\omega}$-ultrafilters is empty.
\end{thm}

\Proof: Suppose that $u$ is a strict $J_{\omega^\omega}$-ultrafilter, thus by definition
of this class, for each $n<\omega$ there exists a function
$f_n:\omega\rightarrow \omega_1$ such that $\ot(f_n(U))\geq
\omega^n$ for each $U \in u$ and there is no function $f_\infty :
\omega \rightarrow \omega_1$ that $\ot(f_\infty(U))\geq
\omega^\omega$ for each $U\in u$. Let $K$ be a set of all such
$n<\omega$ that there is $f_n: \omega \rightarrow \omega_1$ such
that $\ot(f_n(U))\geq \omega^n$ and that there is $U_n\in u$ such
that $\ot(f_n(U_n))=\omega^n$ $^{1)}$\footnotetext{\noindent $^{1)}$ In fact $K =
\omega$, but since we do not need this in the theorem, we omit a
short proof of this fact. }. By Proposition 3.7 there is a sequence
$W_n$ of monotone sequential cascades such that $r(W_{n+1})
> r(W_n)$ and $\int W_n \subset u$.

We will build a sequence $(T_n)$ of monotone sequential cascades
such that

1) $(r(T_n))$ is an increasing sequence

2) $\int T_n \in u$

3) For each $n<\omega$ there exist sets $A_n^0$ and $A_n^1$ that
$A_n^i \# \int T_n$ for $i \in \{0,1\}$ and there is such $i_n \in
\{0,1\}$ that $T_{n+1} \Rrightarrow T_n^{\downarrow A_n^{i_n}}$;

4) $\bigcup_{i\in \{1,\ldots,n\}, j\in\{1,\ldots,i\}} \max T_{i,j}
\cap \max T_{n+1} = \emptyset$, for $T_i = (j) \looparrowleft
T_{i,j}$;

5) $\bigcap_{n<\omega} \max T_n = \emptyset$.

Define $T_1$ as the monotone sequential cascade of rank $1$ with
$\max T_1 = \omega$, clearly $\int T_1 \in u$.  Suppose that
cascades $T_n$ are already  defined for $n\leq m$. For a cascade
$T_{m}= (k) \looparrowleft T_{m,k}$ consider sets $B_{m}^{i} =
\bigcup_{k=2j+i} \max T_{m,k}$ for $i\in \{0,1\}$. Clearly
$B_m^0\#\int T_m$ and $B_m^1 \# \int T_m$, and since $B_m^1 \cup
B_m^0 = \max T_m$ one of these sets belongs to $u$,
call $i_m$ the $i$ for which it happens. Let $C_m = \emptyset$ for $m=2$ and $C_m =
\bigcup_{i\in \{1,\ldots,n\}, j\in\{1,\ldots,i\}} \max T_{i,j}$ for
$m \geq 3$. Clearly $B_m^{i_m} \setminus C_m \in u$. So $\int
T_m^{\downarrow B_m^{i_m} \setminus C_m} \in u$. By Corollary 3.6
applied to $T_m^{\downarrow B_m^{i_m}}$, $W_{m+1}$ and to $u$, there
is a monotone sequential cascade $T_{m+1}$ of finite rank not less then
$m+1$ that $\int T_{m+1} \in u$ and $T_{m+1} \Rrightarrow
(T_m^{\downarrow B_m^{i_m}})$. We define $A_m^{i_m} = \max T_{m+1}$
and $A_m^b= \omega \setminus A_m^{i_m}$ for $b\in \{0,1\}$, $b\not=
i_m$. Clearly $T_{m+1}$ with $A_m^{i_m}$ and $A_m^b$ fulfill the
claim for $n=m+1$. To see that $\bigcap_{n<\omega} \max T_n =
\emptyset$, it suffices to notice that $\bigcup_{k\in  \{1, \ldots,
m\}} \max T_{2,k} \cap \max T_{m+1} = \emptyset$.

Define $T = \bigcup_{n<\omega} T_n^{\downarrow (\omega \setminus
A_{n+1}^{i_{n+1}})}$ ordered by:

1) If $t_1 \in T_n$ and $t_2 \in T_m$  ($n \not= m$) then $t_1
\sqsubset_T t_2 \Leftrightarrow n < m$

2) If $T_1, T_2 \in T_n$ then $t_1 \sqsubseteq_T t_2 \Leftrightarrow
t_1 \sqsubseteq_{V_n} t_2$, where "$\sqsubseteq_{V_n}$" is an order
on $V_n$.

Let $f_\infty: \omega \rightarrow \omega_1$ be a preserving
$\sqsubseteq_T$ order function. Take any $U \in u$ and $n \in
\omega$. Since $ \int T_{n+1} \in u$ thus $U\# \int T_{n+1}$ and so
$U \# \int T_{n+1,k}$ for infinitely many $k$, take $k_0$ from this
set.

Since $U \# \int T_{n+1,k_0}$ thus $\ot(f_{T_{n+1}}(U\cap \max
T_{n+1,k_0})= r(T_{n+1,k_0})\geq \omega^n$. By condition 4) $U\cap
\max T_{n+1,k_0} \cap \max T_i \not= \emptyset $ only for a finite
number of $i>n+1$. So $\{ \max T_{n+1,k_0} \cap U \cap (\max T_i
\setminus \bigcup_{j>i}\max T_j): i\geq n+1\}$ is a finite partition
of $U\cap \max T_{n+1,k_0}$. Thus there is $i_0 \geq n+1$ such that
$\ot\left(f_{T_{n+1}}\left(\max T_{n+1,k_0} \cap U \cap (\max
T_{i_0} \setminus \bigcup_{j>i_0}\max T_j)\right) \right) =$ \\
$\ot\left(f_{T_{n+1}}(\max T_{n+1,k_0} \cap U )\right) \geq
\omega^n$, and since $T_{i_0}\Rrightarrow T_{n+1}$ thus \linebreak
$\ot\left(f_{T_{i_0}}\left(\max T_{n+1,k_0} \cap U \cap (\max T_i
\setminus \bigcup_{j>i_0}\max T_j)\right) \right) \geq \\
\indec \left(\ot\left(f_{T_{n+1}}\left(\max T_{n+1,k_0} \cap U \cap
(\max T_{i_0} \setminus \bigcup_{j>i_0}\max T_j)\right) \right) \right)$,\\
and since $f_\infty\mid_{\max T_{i_0} \setminus \bigcup_{j>i_0}\max
T_j} = f_{T_{i_0}}\mid_{\max T_{i_0} \setminus \bigcup_{j>i_0}\max
T_j}$ thus \\ $\ot\left(f_\infty\left(\max T_{n+1,k_0} \cap U \cap
(\max T_{i_0} \setminus \bigcup_{j>i_0}\max T_j)\right) \right) \geq
\omega^n$.\\ Therefore $\ot(f_\infty(U))\geq \omega^\omega$. $\Box$

\vspace{6mm}
There is a straight correspondence between cascades and $<_\infty$-
sequences.

Let $u$ be an ultrafilter, take sequence $u = u_0 >_\infty u_1
>_\infty \ldots >_\infty u_n$ and functions $f_m: \omega \rightarrow
\omega$ - witnesses that $u_{m-1} >_\infty u_m$.

We will build a monotone sequential cascade $V$ which correspond to
the sequence above with respect to some $U \in u$. In this aim we build a
sequence of cascades $(W_i)_{i\leq n}$. Take any monotone sequential
cascade ${W}_1$ of rank 1 and label elements of $\max {W}_1$ by
natural numbers by any bijections. Clearly  $W \in u_n$ for each $W
\in \int W_1$. Take $W_2 = W_1 \cup \bigcup_{i \in \max W_1,
\card(f_n^{-1}(n))=\infty} f_n^{-1}(i)$ $^{2)}$\footnotetext{\noindent $^{2)}$Since
formally levels in cascade can not intersect we may assume that
domain of $f_1$ and ranges of $f_m$ are subsets of a pairwise disjoint
copies of $\omega$.} ordered by, extended by transitivity, the
following preorder: If $w_1$, $w_2 \in W_1$ then $w_1
\sqsubseteq_{W_2} w_2$ iff $w_1 \sqsubseteq_{W_1} w_2$; if $w_1 \in
\max W_1$ and $w_2 \in \max W_2$ then $w_1 \sqsubseteq_{W_2} w_2$ if
$f_n^{-1}(w_1)=w_2$. Clearly $W \in u_{n-1}$ for each $w \in \int
W_2$. We continue this procedure to get $W_n$ and define $V=W_n$.

Now take any monotone sequential cascade $V$ of finite rank, with
$\int W \subset u$, without loss of generality we may assume that
all branches of $V$ have the same length $n$. For each $v \in V$ let
$\hat{v}$ be an arbitrary element of $\max v^\uparrow$. Consider
functions $f_i:\omega \rightarrow \omega$ such that
$f_1(v_1)=\hat{v}$ for each $v_1 \in \max v^\uparrow $ where
$r(v)=i$. Thus $u >_\infty f_1(u) >_\infty f_2 \circ f_1(u)
>_\infty \ldots >_\infty f_n \circ f_{n-1}\circ \ldots \circ
f_2\circ f_1(u)$, (for details see \cite{Star-Ord-vs-P}).

This cascades - $<_\infty$-sequences correspondence allows us to look
at the Proposition 3.3 and Theorems 3.5 and 3.9 (in virtue of its
proofs) in the following way:

Proposition 3.3 and Theorem 3.5 describes mutual behavior of the
functions - witnesses of $<_\infty$-sequences. Clearly existence of
infinite increasing $<_\infty$-sequences under some ultrafilter
implies existence of an arbitrary long finite $<_\infty$-sequences
under this ultrafilter. Theorem 3.9 shows that if an ultrafilter has
an arbitrary long finite $<_\infty$-sequences then is at least a
strict $J_{\omega^{\omega+1}}$-ultrafilter.

\vspace{3mm}

We'd like to drew attention, not only to benefits, but also to limitations of
the construction presented in the paper. Probably Theorem 3.5 can be proved in a stronger, i.e. infinite version,
but still there is rather no hope to extend our construction to other
limit ordinals. The problem lays in the relations between order ultrafilters and
monotone sequential contours, contained in an ultrafilter,
described in Proposition 3.8, with a special emphasis on
of the upper limitation of the order-type of images. This limitation is non-removable, what was shown
in \cite[Theorem 3.9]{Star-P-hier} by proving (under MA$_{\sigma-centr}$)
that there is a strict $J_{\omega^{\omega+1}}$-ultrafilter that does not contain any monotone
sequential contour of rank 3. Thus we restate Baumgartner question in virtue of our result.

\begin{prob} What about other limit classes? Is there a model with non-empty class of the strict $J_{\omega^{\alpha}}$
ultrafilters for some (all) limit $\omega<\alpha<\omega_1$?  
\end{prob}

Opposite side of this problem is a Shelah question

\begin{prob}\cite[Question 3.12]{Shelah}
Prove the consistency "there is no $J_\alpha$-ultrafilter on $\omega$".
\end{prob}

Under following three theorems of Baumgartner and remembering Shelah model with no P-points,
the above question essentially asks about classes of limit index
and classes whose index is a successor of a limit ordinal.

\begin{thm}\cite[Theorem 4.1]{Baum}
The strict-$J_{\omega^2}$ ultrafilters are P-point ultrafilters.
\end{thm}

\begin{thm}\cite[Theorem 4.2]{Baum}
If there is a P-point then there are strict-$J_{\omega^{\alpha+1}}$ ultrafilters 
for all $\alpha < \omega_1$
\end{thm}

\begin{thm}\cite[Theorem 4.6]{Baum}
Let $\alpha < \omega_1$ and assume $u$ is a strict $J_{\omega^{\alpha+2}}$ ultrafilter.
Then there is a P-point $v$ such that $v \leq_{RK} u$.
\end{thm}

 \hspace{-2mm}
%\gd

%\medskip
%
%{\small\sc \noindent {Wydzia\l}  {Matemetyczno-Fizyczny},
%Politechnika \'{S}l\c{a}ska, Gliwice, Poland
%
%E-mail:  andrzej.starosolski@polsl.pl}

\end{document}